\documentclass[a4paper,12pt]{amsart}

\usepackage{changebar}
\let\cbstart\relax
\let\cbend\relax

\usepackage[latin1]{inputenc}
\usepackage[T1]{fontenc}
\usepackage{url}
\usepackage{subfigure,graphicx}
\usepackage{amsmath,amssymb,mathrsfs}
\usepackage[english]{babel}
\usepackage[all,dvips]{xy}
\usepackage{fleches}

\makeatletter
  \begingroup
    \divide\count0 by 60
    \multiply\count0 by 60
    \advance\count1 by -\count0
    \xdef\dateandtime{%
      \the\year/\two@digits{\the\month}/\two@digits{\the\day}\ 
      \ifnum\count2<10 0\fi \the\count2:%
      \ifnum\count1<10 0\fi \the\count1
    }%
  \endgroup
\makeatother
\usepackage{fancyhdr}
  \fancypagestyle{firstpage}{%
    \fancyhf{}
    \fancyfoot[CO,CE]{\footnotesize\thepage}
    
    }

\makeatletter
\def\@secnumfont{}
\def\section{\@startsection{section}{1}%
  \z@{.7\linespacing\@plus\linespacing}{.5\linespacing}%
  {\normalfont\bfseries\large\centering}}
\makeatother

\usepackage{amsthm}
\theoremstyle{plain}
\newtheorem{theo}{Theorem}[section]
\newtheorem{lemma}[theo]{Lemma}
\newtheorem{coro}[theo]{Corollary}
\newtheorem{prop}[theo]{Proposition}
\theoremstyle{definition}
\newtheorem{rem}[theo]{Remark}
\newtheorem{notations}[theo]{Notations}

\newtheorem{example}[theo]{Example}
\newtheorem{defi}[theo]{Definition}

\def\O{\mathscr O}
\def\CC{\mathbb{C}}
\def\PP{\mathbb{P}}
\def\NN{\mathbb{N}}

\let\phi\varphi
\let\epsilon\varepsilon
\def\M{\mathscr{M}}
\def\C{\mathscr{C}}
\def\Mgn{\M_{g,n}}
\def\Sgn{\mathscr{S}_{g,n}}
\def\Mb{\bar\M}
\def\Mgnb{\bar\Mgn}
\let\Mbgn\Mgnb
\def\gl{_{\text{\textrm{gl}}}}
\def\loc{_{\text{\textrm{loc}}}}

\let\div\relax
\DeclareMathOperator{\div}{div}
\DeclareMathOperator{\sym}{Sym}
\DeclareMathOperator{\spec}{Spec}
\DeclareMathOperator{\Frac}{Frac}
\let\hom\relax
\DeclareMathOperator{\hom}{Hom}
\DeclareMathOperator{\ob}{Ob}\let\ob\relax
\DeclareMathOperator{\pr}{pr}
\DeclareMathOperator{\aut}{Aut}
\DeclareMathOperator{\vect}{Vect}

\def\et{_{\text{\rm ét}}}
\def\ie{\emph{i.e.}}
\let\tens\otimes
\let\wh\widehat
\let\wt\widetilde
\let\leq\leqslant

\let\iso\simeq
\let\equiv\approx
\def\hattens{\mathbin{\widehat{{\otimes}}}}

\makeatletter
\let\OLD@times\times
\def\Times{\@ifnextchar_\timesWSB\OLD@times}
\def\timesWSB_#1{\mathbin{\mathop{\OLD@times}\limits_{#1}}}
\let\OLD@otimes\otimes
\def\Otimes{\@ifnextchar_\otimesWSB\OLD@otimes}
\def\otimesWSB_#1{\mathbin{\mathop{\OLD@otimes}\limits_{#1}}}
\makeatother

\def\d#1{\textbf{\textup{#1}}}

\usepackage{mathrsfs}
\def\fx#1{\mathscr{#1}}

\def\f#1{\mathcal{#1}}

\def\ideal#1{\mathfrak{#1}}

\def\cat#1{\mathfrak{#1}}%

\newcommand{\definefunction}[4]{%
  \ensuremath{
    \left\{
      \begin{array}{ccl} 
        {\displaystyle #1} & \lra & {\displaystyle #2} \\
        {\displaystyle #3} & \lmt & {\displaystyle #4} \\
      \end{array}
    \right.}}

\usepackage{fleches}
\newcommand{\eq}[1][r]
 {\ar@<-3pt>@{-}[#1]
  \ar@<-1pt>@{}[#1]|<{}="gauche"
  \ar@<+0pt>@{}[#1]|-{}="milieu"
  \ar@<+1pt>@{}[#1]|>{}="droite"
  \ar@/^2pt/@{-}"gauche";"milieu"
  \ar@/_2pt/@{-}"milieu";"droite"
  \ar@{}[];[#1]}
\newcommand{\eqv}[1][r]
 {\ar@<1pt>@{}[#1]|<{}="gauche"
  \ar@<2pt>@{}[#1]|-{}="milieu"
  \ar@<3pt>@{}[#1]|>{}="droite"
  \ar@/^2pt/@{-}"gauche";"milieu"
  \ar@/_2pt/@{-}"milieu";"droite"
  \ar@<-3pt>@{}[#1]|<{}="gauche"
  \ar@<-2pt>@{}[#1]|-{}="milieu"
  \ar@<-1pt>@{}[#1]|>{}="droite"
  \ar@/^2pt/@{-}"gauche";"milieu"
  \ar@/_2pt/@{-}"milieu";"droite"
  \ar@{}[];[#1]}
\newcommand{\UN}[4][r]{%
  \ar@/^1pc/[#1]^{#2}_*=<0.3pt>{}="HAUT"
  \ar@/_1pc/[#1]_{#3}^*=<0.3pt>{}="BAS"
  \save\POS "HAUT",*{
      \vrule height 2pt depth 2pt width 0pt
      \vrule height 0pt depth 0pt width 4pt
      }="HAUT",\restore
  \save\POS "BAS",*{
      \vrule height 2pt depth 2pt width 0pt
      \vrule height 0pt depth 0pt width 4pt
      }="BAS",\restore
  \ar @{=>} "HAUT";"BAS" ^{#4}
}
\def\XYMATRIX#1{\raisebox{.5\depth}{\xymatrix{#1}}}
\def\enlargexyentry#1{%
  \POS "#1",*{
    \vrule height 2pt depth 2pt width 0pt
    \vrule height 0pt depth 0pt width 4pt
    }="#1",}
\makeatletter

\def\english{\shorthandoff{;:!?}}
\let\english\relax

\begin{document}
\title[Tangencial base points on algebraic stacks]{Tangencial base points \\ on algebraic stacks}
\author{Vincent Zoonekynd}
\address{Université Paris 7,
UFR de Mathématiques,
Équipe Topologie et Géométrie Algébriques,
Case 7012,
2 place Jussieu,
75251 Paris Cedex 05}
\email{zoonek@math.jussieu.fr}
\urladdr{http://www.math.jussieu.fr/\textasciitilde zoonek/}
\begin{abstract}
  The notion of tangencial base point is well-known for
  schemes in characteristic zero \cite{deligne}. We show
  that the definition in terms of Puiseux series generalizes
  to the case of Deligne--Mumford stacks, over a field of
  arbitrary characteristic. 
  \cbstart
  We then apply this construction to moduli
  spaces of smooth curves, generalizing a result of \cite{IN}.
  \cbend
\end{abstract}
\maketitle

\section{Introduction}

\cbstart
We first outline the construction of a tangencial base point on a
scheme over the field of complex numbers \cite{deligne,AI2}.
\cbend

Let $X$ be a scheme over the field of complex numbers, $D$ a
divisor with normal crossings on~$X$, $x$ a closed point of
$D$ and $\cat{Rev}\,X\setminus D$ the category of étale
coverings of $X\setminus D$. One may ``complete''
$\cat{Rev}\,X\setminus D$ into a topos
$\cat{SRev}\,X\setminus D$, that of disjoint sums of étale
coverings, which is equivalent to the topos of sums of
locally constant étale sheaves on~$X$. A tangencial base
point on $X$ is a functor $\f F : \cat{Rev}\,X\setminus D
\lra \cat{Set}$ whose extension $\f F' : \cat{SRev}\,
X\setminus D \lra \cat{Set}$ defined by $\f F'(\coprod R_i)
= \coprod\f F(R_i)$ is the inverse image functor of a topos
morphism $\cat{Set} \lra \cat{SRev}\,X\setminus D$, \ie, $\f
F'$ is left exact and has a right adjoint.

A choice of coordinates defining $D$ in a formal
neighborhood of $x$, \ie, isomorphisms 
$$\xymatrix{ 
  \spec \wh \O_{X,x} \eq[r] & 
  \CC[\![ t_1, \dots, t_m ]\!] \\
  \spec \wh\O_{D,x} \ar[u] \eq[r] & 
  \dfrac{ \CC[\![ t_1, \dots, t_m ]\!] }{ (t_1\dots t_{m'}) },
  \ar[u] } $$
yields a base point 
$$
\definefunction
{\cat{Rev}\, X \setminus D}
{\cat{Set}}
{R}
{\hom_k(X) \bigl( k(R), \Frac \CC\{\!\{t_1,\dots,t_m\}\!\} \bigr)}
$$
where $\CC\{\!\{t_1,\dots,t_m\}\!\} = 
\limind \CC[\![t_1^{1/n},\dots,t_m^{1/n}]\!]
$
is the ring of Puiseux series 
\cbstart
and $\Frac$ denotes the field of fractions of a ring
\cbend. 

We shall show that this result still holds for schemes over
a field of arbitrary characteristic and generalizes to
Deligne--Mumford stacks. 
\cbstart
We shall then apply it to the special case of
moduli stack of smooth curves. 
The notion of fundamental group of an
algebraic stack will be studied more systematically in a forthcoming
article \cite{Z:pi1}. 
\cbend

\section{Tangencial base points for schemes}

\cbstart
In this section, we prove that the construction outlined for schemes
of characteristic zero  generalizes to schemes over a field of
arbitrary characteristic. 
\cbend

\cbstart
Let us first recall the definition of a divisor with normal crossings from
\cite{GM}. 
\cbend

\begin{defi}
  A divisor $D$ on a scheme $X$ is said to have \d{strictly
    normal crossings} if its support only contains regular
  points of $X$ and if, in the neighborhood of any point
  $s$, its support may be written $D = \sum_{i=1}^{m'} \div
  x_i$, where $(x_i)_{1\leq i \leq m}$ is a regular system
  of parameters at~$s$.

  A divisor $D$ is said to have \d{normal crossings} if, étale
  locally, it has strictly normal crossings, \ie, if there exists a
  surjective étale morphism $U \lre X$ such that 
  $D \times_X U$ be a divisor on $U$ with strictly normal crossings. 
\end{defi}

\begin{defi}
  A \d{Kummer covering} is a morphism
  $$ \spec \dfrac{
    A[t_1, \dots, t_m] 
    }{
    \langle t_1^{n_1} - a_1, \dots, t_m^{n_m} - a_m\rangle
    } \lra \spec A
  $$
  where: $A$ is a ring ; for any prime ideal $\ideal p
  \in \spec A$, the $a_i$ which are invertible in $A_{\ideal
    p}$ are part of a regular system of parameters of
  $A_{\ideal p}$; none of the $n_i$ is divisible by a
  residual characteristic of $A$.
  
  As the $(a_i) _{1\leq i \leq m}$ define a divisor with
  normal crossings, we will say that it is a Kummer covering
  with respect to this divisor.
  
  A \d{tamely ramified covering} over a normal scheme $X$
  with respect to a divisor with normal crossings $D$ is a
  morphism $R \lra X$ that is, étale locally, a disjoint
  sum of Kummer coverings.
  
  We shall denote $\cat{Rev}^D X$ the category of connected
  tamely ramified coverings of $X$ with respect to $D$ and
  $\cat{SRev}^D$ the category (actually, the topos) of their
  disjoint sums.
  In case $D$ is empty, we write 
  $\cat{Rev}\,X$ and $\cat{SRev}\,X$.
\end{defi}

\begin{defi}
  Let $X$ be a normal scheme and $D$ a divisor with normal
  crossings on $X$. A \d{tangencial base point} on
  $X\setminus D$ is a functor $\f F : \cat{Rev}^D X \lra
  \cat{Set}$ whose extension $\f F' : \cat{SRev}^D X \lra
  \cat{Set}$ is the inverse image functor of a topos
  morphism $\cat{Set} \lra \cat{SRev}^D X$.
  In case $D$ is empty, we speak of \d{base point}.
\end{defi}

\begin{rem}
  If $X$ is a normal scheme and $D$ a divisor with normal
  crossings,
  a geometric point $x : \spec K \lra X\setminus D$
  defines a base point on $X\setminus D$:
  $$ \definefunction
  {\cat{Rev}^DX}
  {\cat{Set}}
  {R}
  {\hom_X(\spec K, R).}
  $$
\end{rem}

\begin{defi}
  Let $A$ be a regular local (hence integral) ring, $\ideal
  m$ its maximal ideal, $(t_1,\dots, t_m)$ a regular system
  of parameters, $\wh A = \limpro A/\ideal m^n$ its $\ideal
  m$-adic completion and $\bar \Omega$ an algebraic closure
  of $\Frac \wh A$.  For each $1\leq i \leq m$ and each $n$
  prime to the characteristic of $\bar \Omega$, we choose a
  primitive $n$th root $t_i^{1/n}$ of $t_i$ in $\bar
  \Omega$, and we define the \d{Puiseux ring} of $A$ as
  $$ \wt A_{t_1, \dots, t_m} = 
  \limind_n \wh A[t_1^{1/n}, \dots, t_m^{1/n}]$$
  where the morphisms defining this inductive system are 
  $$\definefunction
  {\wh A[t_1^{1/k}, \dots, t_m^{1/k}]}
  {\wh A[t_1^{1/k\ell}, \dots, t_m^{1/k\ell}]}
  {t_i^{1/k}}
  {(t_i^{1/k\ell})^\ell.}
  $$  
  We shall often write $\wt A$ instead of 
  $\wt A_{(t_1,\dots,t_m)}$.
\end{defi}
\begin{rem}
  As a ring, $\wt A$ is merely the sub-$\wh A$-algebra of $\bar
  \Omega$ generated by the $t_i^{1/n}$. But we shall regard
  it as a ring endowed with an inclusion $A \lrh \wt
  A$ and a \d{coherent system} of roots of the $t_i$,
  denoted $t_i^{1/n}$: it means that for all $k$, $\ell$
  prime to the characteristic of $A$, one has
  $(t_i^{1/k\ell})^\ell = t_i^{1/k}$.  These data $(\wt A, A
  \lra \wt A, (t_i^{1/n})_{i,n})$ are well determined, from $A$
  and $(t_1,\dots,t_m)$, up to a
  \emph{unique} isomorphism.
\end{rem}

\begin{theo}
  \label{theo:basepointschemes}
  Let $X$ be an integral scheme, $D$ a divisor with normal
  crossings on~$X$, $x \in D$ a closed point whose residue
  field is algebraically closed and $(t_1, \dots, t_m)$ a
  regular system of parameters of $\O_{X, x}$ such that $D$
  be defined around $x$ by the equation $t_1 \cdots t_{m'} =
  0$.  Then, the functor
  $$ \f F : \definefunction
  {\cat{Rev}^D X}
  {\cat{Set}}
  {R}
  {\hom_{k(X)} \bigl( k(R), \Frac \wt\O_{X, x, (t_1,\dots,t_m)} \bigr)}
  $$
  is a base point.
\end{theo}

\begin{proof}
  Let $\Omega=\Frac \wt\O_{X,x}$ and $\bar \Omega$ an
  algebraic closure of $\Omega$. We shall show that the
  functor $\f F$ is isomorphic to that defined by the
  geometric point $\spec \bar \Omega \lra X \setminus D$.
  Let $R \in \ob \cat{Rev}^D X$ be a tamely ramified
  covering of~$X$ and set $\spec B = R \times_X \spec
  \O_{X,x}$. Let us show the following isomorphisms.
  
  \begin{align}
    \hom_X(\spec \bar\Omega, R) 
    &\iso \hom_{\spec \O_{X,x}}( \spec\bar\Omega, \spec B) 
    \label{eq:1}\\
    &\iso \hom_{\spec \O_{X,x}}( \spec\wh\O_{X,x}, \spec B) 
    \label{eq:2}\\
    &\iso \hom_{\O_{X,x}}(B, \wh\O_{X,x})
    \label{eq:3}\\
    &\iso \hom_{\Frac \O_{X,x}}(\Frac B, \Omega)
    \label{eq:4}\\
    &\iso \hom_{k(X)} \bigl(k(R), \Omega\bigr)
    \label{eq:5}
  \end{align}
  
  The first isomorphism \eqref{eq:1} merely states the
  universal property of the fiber product
  $\spec B = R \times_X \spec \O_{X,x}$.
  $$ \xymatrix{\relax
    &&&\spec B\ar[d]\ar[r]\cartesien&R\ar[d]\\
    \spec\bar\Omega \ar[r] \ar@{.>}[rrru]\ar@/^2.5pc/[rrrru]&
    \spec\Omega \ar[r] &
    \spec \wh\O_{X,x} \ar[r] &
    \spec \O_{X,x} \ar[r] &
    X} $$
  
  To prove \eqref{eq:2},
  set $\wh B = \spec B \times_{\spec \wh \O_{X,x}} 
  \spec \O_{X,x}$.
  As $R \lra X$ is a tamely ramified covering, 
  $\spec \wh B \lra \spec \wh \O_{X,x}$ is a disjoint sum of
  Kummer coverings. 
  The morphism 
  $\spec \bar \Omega \lra \spec \wh B$ determines a
  connected component 
  $\spec \wh B _i$ of $\spec \wh B$ and the image of 
  $\wh B_i \lra \bar \Omega$ is contained in 
  $\wt \O_{X,x}$, hence the morphism 
  $\spec \bar \Omega \lra \spec \wh B$ extends to 
  $\spec \wt \O_{X,x} \lra \spec \wh B$. 
  $$ \xymatrix{
    && \spec \wh B_i \ar[d] \ar[r] 
    & \spec \wh B \ar[d] \\
    \spec \bar\Omega \ar[r] \ar[rru] &
    \spec \wt \O_{X,x} \ar@{.>}[ru] \ar[r] &
    \spec \wh\O_{X,x} \ar[r] & 
    \spec \wh\O_{X,x}.
    }
  $$

  Isomorphisms \eqref{eq:3} and 
  \eqref{eq:5} are straightforward, and 
  \eqref{eq:4} results from lemma \ref{lemma:point4}.
\end{proof}

\begin{lemma}
  \label{lemma:proptild}
  Let $A$ be a noetherian, local (hence of finite
  dimension), regular (hence integral) ring with
  algebraically closed residue field $k$, 
  let $(t_1,\dots,t_m)$ be a regular system of parameters of
  $A$ and let $B$ be a
  finite integral $A$-algebra. 
  Then $\wt B \iso \wt A \tens _A B$ and $\wt B$ is a finite
  $\wt A$-algebra.
\end{lemma}
\begin{proof}
  One has 
  \begin{align*}
    \wt B &=
    \limind \wh B [x_1^{1/m},\dots,x_n^{1/m}] \\
    &\iso \limind (B\otimes_A \wh A)
    [x_1^{1/m},\dots,x_n^{1/m}] \\ 
    &\iso\limind B\otimes_A 
    \wh A[x_1^{1/m},\dots,x_n^{1/m}] \\ 
    &\iso B\otimes_A \limind \wh
    A[x_1^{1/m},\dots,x_n^{1/m}] 
    &&\text{from \cite[7.2]{E}}\\ 
    &\iso B\otimes_A \wt A.
  \end{align*}
\end{proof}
\begin{lemma}
  \label{lemma:point4}
  Let $A$ be a noetherian, local (hence of finite
  dimension), regular (hence integral) ring with
  algebraically closed residue field $k$, 
  let $(t_1,\dots,t_m)$ be a regular system of parameters of
  $A$ and let $B$ be a
  finite integral $A$-algebra. 
  Then 
  $\hom_A(B,\wt A) \iso 
  \hom_{\Frac A}(\Frac B, \Frac \wt A)$.
\end{lemma}
\begin{proof}
  We shall show, by induction on the dimension of $A$, that
  all $A$-morphisms $\phi : B \lra \wt A$ are injective. 

  The case $\dim A = 0$ is trivial, for then $A = B = k$.


  Assume the result holds for rings of dimension $n-1$ and
  let $A$ be of dimension $n$. 
  Let $\phi : B \lra \wt A$ be any $A$-morphism and set 
  $\ideal p = \ker \phi$. 
  As $B/\ideal p$ is isomorphic to a subring of $\wt A$, it
  is integral and thus $\ideal p$ is prime. 
  From \cite[10.14]{E}, there exists $t \in A$ such that
  $(t)$ be prime and $A/(t)$ be a local regular ring of
  dimension $n-1$.  One may then write 
  $$\xymatrix{\relax \smash[b]{B \Otimes_A \bigl(A/(t)\bigr)
      \iso B/tB}\vphantom{B} \ar[rr]^{\phi'} &&
    \wt A \Otimes_A A/(t) \iso \wt{A/(t)} \\
    & A/(t). \ar[ru]\ar[lu]}$$
  where the isomorphism on the right comes from lemma 
  \ref{lemma:proptild}. By the induction assumption, $\ker \phi'
  = 0$, hence $\ideal p + tB = tB$, hence 
  $\ideal p \subset tB$.
  In geometric terms, $\spec B$ is a scheme, 
  $\spec B/tB$ is a hypersurface and $\spec B/\ideal p$ is
  an irreducible subscheme containing this hypersurface: it
  is either all of $\spec B$ or $\spec B/tB$. As $t$ is
  nonzero in $B/\ideal p$, we must reject the first
  possibility: $\spec B/\ideal p = \spec B$, hence $\ideal p
  = 0$.
\end{proof}

\section{Tangencial base points for stacks}

\cbstart
In this section, we generalize to algebraic stacks the preceeding
construction of tangencial base points. 
\cbend
By \d{stack} we shall always mean Deligne--Mumford stack:
see \cite{DM}, \cite{V}, \cite{LMB} or \cite{Z}.

\begin{rem}
  The \d{étale site} $X\et$ of a stack~$X$ is defined in
  \cite{DM}. We shall denote $\cat{Sh}\,X$ the corresponding
  category of sheaves.
  Let $\cat{Rev}\,X$ be the
  category of connected locally constant finite étale
  sheaves on~$X$ and $\cat{SRev}\,X$ the category of their
  disjoint sums (it is a topos). This category could also be
  defined as a category whose objects are finite étale stack
  morphisms $R \lra X$, but, because stacks form a $2$-category,
  the morphisms are more
  complicated to define: see \cite{Z:ThVK}.
  
  If $X$ is a stack, there exists a surjective étale
  morphism $U \lre X$ from a scheme $U$ and one may prove
  (see \cite{V}) that $U \times_X U \dlra U$ is an étale
  groupoid (\ie, a groupoid whose source morphism is étale)
  whose diagonal $U \times_X U \lra U \times U$ is quasi
  compact and separated. Conversely, any such groupoid
  defines a quotient stack, denoted $[U/S]$. 

  Let $\cat{Rev}(S\dlra U)$ denote the category of connected
  equivariant étale coverings: its objects are étale
  coverings $R \lra U$ endowed with a morphism 
  $\alpha : 
  S 
  \mathbin{{}_s{\Times_U}}
  R 
  \lra R$ such that 
  $$\xymatrix{
    S \times_V R
    \ar@<+.3ex>[r]^-{\pr_2}
    \ar@<-.3ex>[r]_-{\alpha}
    \ar[d] &
    R \ar[d] \\
    S \dar[r] & U } $$
  be a groupoid morphism and such that the (set-theoretic)
  groupoid $\pi_0( R \times_V S) \dlra \pi_0 R$ be
  connected; its morphisms 
  $(R,\alpha) \lra (R',\alpha')$ are the morphisms 
  $f : R \lra R'$ such that the following diagrams commute. 
  $$ \xymatrix{
    R \ar[rr]^f \ar[rd] && R' \ar[ld] \\ & V
    } 
  \qquad\qquad
  \xymatrix{ 
    R \times _V S \ar[r]^-\alpha \ar[d]_{f \times 1} & 
    R \ar[d]^f \\ 
    R' \times_V S \ar[r]_-{\alpha'} & R' }
  $$
  Similarily, we let 
  $\cat{SRev}(S\dlra U)$ denote the category of disjoint
  sums of equivariant étale coverings (it is a topos). 
 
  One may show (see \cite{V}) that if $(S \dlra U)$ is a
  presentation of the stack~$X$, then there is an
  equivalence of categories 
  $$\cat{Rev}\,X \equiv \cat{Rev}(S\dlra U). $$
\end{rem}

\begin{lemma}
  A stack morphism $f: Y \lra X$ induces topos morphisms
  $\cat{Sh}\,Y \lra \cat{Sh}\,X$ 
  and 
  $\cat{SRev}\,Y \lra \cat{SRev}\,X$.
\end{lemma}
\begin{proof}
  As the functor 
  $$ f^* : \definefunction
  {X\et}{Y\et}
  {\raisebox{.5\depth}{\xymatrix{U\ar[d]\\X}}}
  {\raisebox{.5\depth}{\xymatrix{U \times_X Y \ar[d] \\ Y }}}
  $$
  is continuous, \ie, as it 
  preserves finite projective limits 
  and covering sieves (if $U \lra X$ is an object of 
  $X\et$ and $R = (U_i \lra U)_{i\in I}$ a covering 
  sieve, then the sieve (generated by) 
  $f^* R = (U_i \times_X Y \lra U \times _X Y)$
  is a covering sieve of $U \times_X  Y \lra Y$), 
  it defines by \cite[VII.10.2]{MM:SGL} a topos morphism. 

  The functor $f^*$ clearly restricts to the categories of
  locally constant sheaves
  $\cat{SRev}\,Y\lra\cat{SRev}\,X$ and defines a topos morphism.
\end{proof}

\begin{defi}
  A stack morphism $D \lra X$ is said to be a \d{divisor
    with normal crossings} if there exists a surjective
  étale morphism $U \lra X$ from a scheme such that $D
  \times_X U \lra U$ be a divisor with normal crossings.
  This property then holds for any surjective étale morphism
  from a scheme.
\end{defi}

\begin{defi}
  Let $X$ be a stack and $D$ a divisor with normal crossings
  on $X$. An étale sheaf~$\fx F$ on $X$ is said to be a
  \d{tamely ramified covering} of~$X$ along $D$
  if there is a surjective étale morphism 
  $\pi: U\lra X$ from a
  scheme $U$ such that $\pi^* \fx F$ be a tamely ramified
  covering of~$U$ along $D \times_X U$. We denote by 
  $\cat{Rev}^D X$ the category of connected tamely ramified
  coverings of~$X$ along $D$ and $\cat{SRev}^DX$ the category
  of their disjoint sums (it is a topos). 

  A stack morphism $Y \lra X$ induces a topos morphism
  $ \cat{SRev}^{D\times_XY} Y \lra 
  \cat{SRev}^DX$.

  A \d{tangencial base point} on $X\setminus D$ is a functor
  $\cat{Rev}^DX\lra \cat{Set}$ whose extension to 
  $\cat{SRev}^DX\lra \cat{Set}$ is the inverse image functor
  of a topos morphism.
\end{defi}

\begin{rem}
   If $D$ is a divisor with normal crossings on a stack~$X$, we denote
   $\cat{Rev}^D(S\dlra U)$ the category of connected 
   equivariant tamely ramified coverings of~$U$ and 
   $\cat{SRev}^D(S\dlra U)$ the category of their disjoint
   sums.  
   As before, one has an equivalence of categories
   $\cat{Rev}^D\,X \equiv \cat{Rev}^D(S\dlra U)$. 
\end{rem}

\begin{lemma}
  \label{lemma:basepointeasystack}
  Let $X$ be a connected algebraic stack, 
  and $S \dlra U$ a presentation of~$X$. 
  Then, a geometric point 
  $y : \spec K \lra U$ defines a tangential base point on $X$
  $$\f F : \definefunction
  {\cat{Rev}^D(S \dlra U)}
  {\cat{Set}}
  {R}
  {\hom_U(\spec K, R).}
  $$  
\end{lemma}
\begin{proof}
  We have stack morphisms
  $$ \spec K \lra U \lra X, $$
  hence topos morphisms
  $$ \cat{Set} \lra \cat{SRev}^{D\times_XU}U \lra 
  \cat{SRev}^D X \equiv \cat{SRev}^D(S \times U \dlra U) 
  $$
  where $S\times U \dlra U$ is a presentation de $X$. 
  The inverse image of the second is the forgetful functor
  (forgetting the action of the groupoid) and that of the
  first is  
  $$\definefunction
  {\cat{SRev}\,\spec\wh\O_{U,y}}
  {\cat{Set}}
  {R}
  {\hom_{\spec\wh\O_{U,y}}(\spec K, R).}
  $$
  The composition of these functors yields $\f F$.
\end{proof}
\begin{theo}
  Let $X$ be a regular stack over an
  algebraically closed field~$k$, 
  let $S \dlra U$ be a presentation of~$X$, 
  let $D \subset X$ be a divisor with normal crossings, 
  let $y : \spec k \lra D \times_XU$ be a point of 
  $D\times_XU \subset U$, 
  let $x : \spec k \lra U \lra X$ be the corresponding point
  of~$X$, 
  let $(t_1,\dots,t_m)$ be a regular system of parameters
  for $D\times_XU$ at $y$.
  Then the functor 
  $$ \definefunction 
  {\cat{Rev}^D\,X}
  {\cat{Set}}
  {R}
  {\coprod_i \hom_{ k(U_0) }
    \bigl( 
    k(R_i), 
    \Frac \wt \O_{U,y, (t_1,\dots,t_m)}
    \bigr)}
  $$
  is a base point, where 
  $U_0$ is the connected component of~$U$ containing $y$ and
  the $R_i$ are the connected components of 
  $R \times_X U_0$.  
\end{theo}
\begin{proof}
  Let $R \in \ob\cat{Rev}^D X$ be a tamely ramified covering
  of~$X$ and $\bar \Omega$ an algebraic closure of 
  $\Frac \wh \O_{U,y}$. From lemma
  \ref{lemma:basepointeasystack},  
  we know that the functor 
  $$ \definefunction
  {\cat{Rev}^D X}
  {\cat{Set}}
  {R}
  {R \times_X U}
  $$
  is a base point. But 
  \begin{align*}
    \hom_U\bigl(\spec \bar \Omega, R\times_X U\bigr) 
    &\iso \hom_{U_0}\bigl(\spec \bar \Omega, R\times_X U_0\bigr) \\
    &\iso \hom_{U_0}\bigl(\spec \bar \Omega, \coprod R_i\bigr) \\
    &\iso \coprod \hom_{U_0}\bigl(\spec \bar \Omega, R_i\bigr) \\
    &\iso \coprod \hom_{k(U_0)}\bigl(k(R_i), \bar \Omega\bigr),
  \end{align*}
  the last isomorphism coming from theorem
  \ref{theo:basepointschemes}. 
\end{proof}

\section{Tangent space of a stack}

\cbstart
The tangencial base points of the introduction may be defined by the
datum of a tangent vector which is not tangent to the divisor or,
equivalently, by the datum of coordinates $t_1,\dots,t_n$ on the
tangent space such that the trace of the divisor be defined by the
equation $t_1 \cdots t_m = 0$. In this section, we define the notion
of tangent space of an algebraic stack and generalize this remark to
the case of stacks. 
\cbend

\begin{lemma}
  Let $X$ be a stack over an algebraically closed field~$k$
  and $x: \spec k \lra X$ a point of~$X$. Then, there is a
  finite group $G$ such that the square 
  $$ \xymatrix{
    G \times \spec k \ar[r]\ar[d]\cartesien &
    \spec k \ar[d]^x \\
    \spec k \ar[r]_x & X. } $$
  be $2$-cartesian.
  We say that $G$ is the \d{fundamental group} of the point~$x$.
\end{lemma}
\begin{proof}
  The existence of the set $G$ is left to the reader:
  proceed as in proposition \ref{prop:localstrstack}. 
  Further, $G \times \spec k \dlra \spec k$ is a groupoid
  (and hence $G$ is a group) for the same reasons that a
  surjective étale morphism from a scheme 
  $U \lra X$ yields a groupoid 
  $U \times_X U \dlra U$: see \cite{V}.
\end{proof}

\begin{prop}
  \label{prop:localstrstack}
  Let $X$ be a stack on an algebraically closed field~$k$, 
  let $U \lre X$ be a surjective étale morphism from a
  scheme, 
  let $y : \spec k \lra U$ be a point of $U$, 
  let $x : \spec k \lra U \lra X$ be the corresponding point
  of~$X$
  and let $G$ be the automorphism group of $x$.
  Then, one has a $2$-cartesian square 
  $$ \xymatrix{
    G \times \spec \wh \O_{U, y} \ar[r] \ar[d] \cartesien 
    & \spec \wh \O_{U,y} \ar[d] \\
    \spec \wh \O_{U,y} \ar[r] & X. } $$
\end{prop}
\begin{proof}
  In the following diagram, the bold arrows are étale. 
  $$ \xymatrix{
    \spec\wh\O_{U,y} \Times_U R \Times_U \spec\wh\O_{U,y} 
    \cartesien
    \ar[r] \ar[d] &
    \smash{R \Times_U \spec\wh\O_{U,y}}\vphantom{\spec}
    \ar[r] 
    \ar@<.4pt>[r]
    \ar@<-.4pt>[r]
    \ar[d] \cartesien &
    \spec\wh\O_{U,y} \ar[d] \\
    \smash[b]{\spec\wh\O_{U,y} \Times_U R}\vphantom{\spec}
    \ar[r] 
    \ar[d] 
    \ar@<.4pt>[d]
    \ar@<-.4pt>[d]
    \cartesien & 
    R 
    \ar[d] 
    \ar@<.4pt>[d]
    \ar@<-.4pt>[d]
    \ar[r] 
    \ar@<.4pt>[r]
    \ar@<-.4pt>[r]
    \cartesien &
    U 
    \ar[d]
    \ar@<.4pt>[d]
    \ar@<-.4pt>[d]
    \\
    \spec\wh\O_{U,y} 
    \ar[r]
    &
    U 
    \ar[r]
    \ar@<.4pt>[r]
    \ar@<-.4pt>[r]
    & 
    X } $$
  Hence, as $\wh\O_{U,y}$ is a strict henselian ring, there
  exists sets $A$ and $B$ such that the diagram reads
  $$ \xymatrix{
    \spec\wh\O_{U,y} \Times_U R \Times_U \spec\wh\O_{U,y} 
    \cartesien
    \ar[r] \ar[d] 
    \ar@{}[rd]|(.65)*{(1)}&
    A \times \spec\wh\O_{U,y} 
    \ar[r] 
    \ar[d] \cartesien &
    \spec\wh\O_{U,y} \ar[d] \\
    B \times \spec\wh\O_{U,y}
    \ar[r] 
    \ar[d] 
    \cartesien & 
    R 
    \ar[d] 
    \ar[r] 
    \cartesien &
    U 
    \ar[d]
    \\
    \spec\wh\O_{U,y} 
    \ar[r]
    &
    U 
    \ar[r]
    & 
    X } $$
  The square $(1)$ may be written, replacing $R$ by 
  the union of the `images' of 
  $A \times \spec \wh\O_{U,y}$
  and
  $B \times \spec \wh\O_{U,y}$, 
  $$ \xymatrix{
    \spec\wh\O_{U,y} \Times_U R \Times_U \spec\wh\O_{U,y} 
    \ar[d] \ar[r] 
    \cartesien
    &
    A \times \spec \wh\O_{U,y}
    \ar[d] \\
    B \times \spec \wh\O_{U,y} 
    \ar[r] & 
    \displaystyle\coprod_{r \in D}
    \spec \wh \O_{R,r}
    }
  $$
  where $D \subset R(k)$. 
  Hence 
  $\spec\wh\O_{U,y} \Times_U R \Times_U \spec\wh\O_{U,y}
  \equiv C \times \spec \wh \O_{U,y}$ for some set $C$. 
  But as 
  $$\xymatrix{
    G \times \spec k \ar[r] \ar[d] \cartesien & 
    C \times \spec k \ar[d] \ar[r] \cartesien &
    \spec k \ar[d] \\
    C \times \spec k \ar[r] \ar[d] \cartesien & 
    C \times \spec\wh\O_{U,y} \ar[r] \ar[d] \cartesien &
    \spec\wh\O_{U,y} \ar[d] \\
    \spec k \ar[r] & \spec\wh\O_{U,y} \ar[r] & X,
    } $$
  where the large square is that defining $G$, one has $G = C$. 
\end{proof}
\begin{defi}
  Let $X$ be a stack on an algebraically closed field $k$, 
  let $x : \spec k \lra X$ be a point of~$X$
  and $G$ its group, 
  let $U \lra X$ be a surjective étale morphism from a
  scheme, 
  let $y : \spec k \lra U$ be a point of $U$ above $x$.
  The \d{tangent space} of~$X$ at~$x$, denoted 
  $T_xX$, is the tangent space 
  $T_y U$ endowed with the action of~$G$. 
\end{defi}
\begin{rem}
  \cbstart
  The tangent space of a stack is not unique up to a unique
  isomorphism: this should not be surprising, for it satisfies the
  following \d{2-universal} property. 
  First, one may identify a
  $k$-vector space~$V$ endowed with an action if a finite group $G$
  with the quotient stack 
  $( \spec \wh{\sym V} )/G$. 
  The tangent space $T_x X$ is then a 2-final diagram among diagrams
  of the form 
  $$ \xymatrix @! C=1cm {
    & \spec k \ar[rd]^x |*{}="A" 
    \ar[ld] 
    \enlargexyentry A
    \\ 
    \spec \wh{\sym V}/G 
    \ar@{=>}"A" \ar[rr]_f && X }
  $$
  where $V$ is a (non-zero) $k$-vector space and $G$ a finite group
  acting linearly on $V$. 
  \cbend
\end{rem}
\begin{rem}
  This notion of tangent space of a stack is equivalent to
  that defined in \cite{LMB} as the category 
  $$
  \cat{hom}^{\mathstrut}_{\spec k}
  \left( 
    \spec \dfrac{k[\epsilon]}{\langle\epsilon^2\rangle}
    ,\ 
    X
  \right)
  $$
  whose objects are the diagrams 
  $$ \xymatrix{
    & \spec k \ar[ld] \ar[rd]^x \\
    \spec \dfrac{k[\epsilon]}{\langle\epsilon^2\rangle}
    \ar[rr]_f && X }
  $$
  and whose morphisms $f \lra g$ are stack $2$-morphisms 
  $$ \xymatrix{
    *!!<0pt,0.7ex>+!R!<-1.2ex,0pt>
    {\spec \dfrac{k[\epsilon]}{\langle\epsilon^2\rangle}}
    \UN{f}{g}{} & 
    *!!<0pt,0.7ex>+!L!<+1.2ex,0pt>
    {X.} 
    }
  $$
\end{rem}
\begin{example}\label{ex:tangentspacetangbasepoint}%
  Let $X$ be an algebraic stack over an algebraically closed
  field $k$, let $D$ be a divisor with normal crossings on
  $X$ and let 
  $x : \spec k \lra D$ be a point of $D$. Coordinates on the
  tangent space $T_x X$ such that the trace of $D$ on 
  $T_xX$ be given by the equation
  $t_1\cdots t_{m'}=0$ define a tangencial base point on $X
  \setminus D$, thereby generalizing the classical situation
  presented in the introduction. 
\end{example}


\cbstart
\section{Tangencial base points on the moduli stack of smooth curves}

We shall now apply the preceeding theory to the case of moduli stacks
of smooth stable pointed curves. 

\begin{defi}
  A \d{stable $n$-pointed curve} of genus $g$ 
  is a proper and flat morphism 
  $\pi : C \sra S$ 
  endowed with $n$ sections
  $s_i : S \sra C$ 
  (for $1\leq i\leq n$) such that 
  the geometric fibers~$C_s$ of $\pi$ 
  over any geometric point $s$ of $S$ 
  are reduced, connected
  curves, of arithmetic genus 
  $\dim H^1(C_s, \O_{C_s}) = g$, whose singular points are ordinary
  double points;
  such that the image of the sections $s_i$ in the geometric
  fibers~$C_s$ are smooth distinct points $P_i$;
  and such that these geometric fibers have a finite number of
  automorphisms fixing the $P_i$. 

  A \d{special point} of a stable pointed curve $C$ is a point of the
  normalization $\wt C$ above a singular or marked point of $C$. 
  A \d{maximally degenerate stable pointed curve} is a stable pointed
  curve $C$ such that the irreducible components of the normalisation
  $\wt C$ be projective lines with three special points. 
\end{defi}

\begin{notations}
  Let us denote $\Mgnb$ the moduli stack of stable $n$-pointed genus
  $g$ curves, $\Sgn$ the divisor with normal crossings of singular
  curves and $\Mgn = \Mgnb \setminus \Sgn$ the moduli stack of smooth
  stable $n$-pointed genus $g$ curves \cite{knudsen}.
\end{notations}

We shall first recall the theory of deformations of curves
\cite{schlessinger,DM}; then explain how they lead to a description of
the tangent space to the moduli stack of stable pointed curves
$\Mgnb$; and finally describe tangencial base points on the moduli
stack of smooth curves $\Mgn$.

We shall state the results above a 
field of characteristic zero: in positive characteristic, it suffices
to replace ``local Artin $k$-algebra with residue field $k$'' by
``local Artin $A_k$-algebra with residue field $k$, where $A_k$ is a
complete regular local ring, with maximal ideal $p A_k$ and residue
field $k$'' \cite[p.~79]{DM}.

\begin{defi}%
  A \d{deformation} of a proper curve 
  $C \lra \spec k$ is a cartesian square 
  $$ \xymatrix{ 
    C \cartesien \ar[r] \ar[d] & \C \ar[d] \\ \spec k \ar[r] & \spec
    A} $$
  where $\C \lra \spec A$ is a (proper and flat) curve and 
  $A$ a local Artin $k$-algebra with residue field~$k$. 

  A \d{deformation morphism} 
  $$ \XYMATRIX{C \ar[r]\ar[d]\cartesien & \C \ar[d] \\ 
    \spec k \ar[r] & \M } 
  \lra 
  \XYMATRIX{C \ar[r]\ar[d]\cartesien & \C' \ar[d] \\ 
    \spec k \ar[r] & \M' } 
  $$
  is a commutative diagram 
  $$ \english 
  \xymatrix@R=5pt{ & C \ar[rd] \ar[ld] \ar[dd]|!{[ld];[dr]}\hole \\
    \C \ar[rr] \ar[dd] && \C' \ar[dd] \\
    & \spec k \ar[rd] \ar[ld] \\
    \M \ar[rr] && \M'. } $$
  Hence, deformation \emph{iso}morphisms are morphisms of families of
  curves whose restriction to $C$ is the identity. 

  A \d{generalized deformation} of a proper curve 
  $C \lra \spec k$ is a cartesian square
  $$ \xymatrix{ 
    C \cartesien \ar[r] \ar[d] & \C \ar[d] \\ \spec k \ar[r] & \spec
    A} $$
  where $\C \lra \spec A$ is a (proper and flat) curve and $A$ a local
  $k$-algebra with residue field $k$ and maximal ideal $\mathfrak m$,
  such that all the quotients 
  $A / \mathfrak m ^n$, with $n\in \NN^\times$, be Artin
  $k$-algebras. Morphisms of generalized deformations are defined as 
  above. 

  A \d{universal deformation} of a curve 
  $C \lra \spec k$ is a generalized deformation 
  $$ \xymatrix{ C \ar[r] \ar[d] & \C \ar[d] \\ \spec k \ar[r] & \M }
  $$
  such that for any deformation 
  $$ \xymatrix{ C \ar[r] \ar[d] & \C' \ar[d] \\ \spec k \ar[r] & \M', }
  $$
  there exists a unique morphism 
  $f : \M' \lra \M$ and a unique morphism 
  $\wt f : \C' \lra \C$ 
  such that the following diagram commute and the front square be
  cartesian. 
  $$ \english \xymatrix@R=5pt{\relax
    & C \ar[rd]\ar[ld]\ar[dd]|!{[ld];[dr]}\hole  \\
    \C' \ar@{.>}[rr] \ar[dd] \cartesien && \C \ar[dd] \\
    & \spec k \ar[rd] \ar[ld] \\
    \M' \ar@{.>}[rr]_f && \M} $$
  The scheme $\M$ is then called the \d{base} of the universal deformation. 
  This is equivalent to requesting that the functor 
  $$ D : \definefunction{\cat{Art}\,k}{\cat{Ens}}A%
  {\{\,\text{isomorphism classes of deformations of $C$ over $A$}\,\}}$$
  be \d{prorepresentable}, where $\cat{Art}\,k$ denotes the category
  of local Artin $k$-algebras with residue field $k$ and their local
  morphisms. 
\end{defi}

\begin{lemma}
  Stable curves have universal deformations. 
\end{lemma}
\begin{proof}
  See \cite[p.~80]{DM}.
\end{proof}

\begin{defi}
  Let $C$ be a stable curve over $k$ and $x$ a singular point of~$C$.
  A \d{local deformation} of $C$ at $x$ over 
  an Artin local $k$-algebra with residue field $k$ is the datum of a
  flat $A$-algebra $\O$ and an isomorphism
  $\O \tens_A k \iso \wh \O_{X,x}$, \ie, 
  the datum of a cartesian square 
  $$ \xymatrix{\relax
    \spec \wh\O_{X,x} \ar[r]\ar[d]\cartesien & \spec \O \ar[d] \\ 
    \spec k \ar[r] & \spec A} $$
  where the morphism 
  $\spec\O \lra \spec A$ is flat.
    
  One says that $C$ has a \d{universal local deformation} at $x$ if the
  functor 
  $$
  D\loc : 
  \definefunction {\cat{Art}\,k}{\cat{Ens}}A%
  {
    \vtop{\hsize=8cm\raggedright\noindent
      \hbox{$\{\,$isomorphism classes of local}
      \hbox{\hphantom{$\{\,$}deformations of $C$ at $x$ over $A \,\}$}}
    }
  $$
  is prorepresentable; 
  a ring (or its spectrum) prorepresenting this fuctor 
  is said to be the \d{base} of the
  universal local deformation of $C$ at $x$. 
\end{defi}

\begin{rem}
  Let $C$ be a stable curve, 
  $x_1,\dots,x_n$ its singular points and, for all $1\leq i\leq n$, 
  $A_i$ the base of the universal local deformation of $C$ at
  $x_i$ \cite[p.~81]{DM}: if $(u_i, v_i)$ are local coordinates
  around~$x_i$, \ie, if 
  one chooses an isomorphism
  $$ \wh \O_{C,x_i} \iso \dfrac{k[\![u_i,v_i]\!]}{u_iv_i}, $$
  the deformation is given by
  $$\xymatrix{\relax
    \spec \dfrac{k[\![u_i,v_i]\!]}{\langle u_iv_i\rangle}
    \ar[r]\ar[d]\cartesien &
    \spec \dfrac{k[\![u_i,v_i,\epsilon]\!]}{\langle u_iv_i-\epsilon_i\rangle}
    \ar[d] \\
    \spec k \ar[r] & \spec k[\![\epsilon_i]\!]
    \rlap{\ensuremath{{}=\spec A_i.}} } $$  
  If $C$ is a maximally degenerate curve,
  the \d{base} of the universal local deformation of $C$ is defined as
  $$\M\loc = \spec M\loc = \spec A_1 \hattens \cdots \hattens A_n
  = \spec k[\![\epsilon_1,\dots,\epsilon_n]\!].$$ 
\end{rem}

\begin{lemma}\label{lemme:defgl=defloc}%
  If $C$ is a maximally degenerate stable curve,
  the base of its universal deformation 
  $\M\gl = \spec M\gl$ is canonically isomorphic
  to the base of its universal local deformation
  $\M\loc = \spec M\loc$.
\end{lemma}
\begin{proof}
  See \cite[1.5 p.80--81]{DM}.
\end{proof}

\begin{lemma}%
  \label{lemme:defunivauto}%
  If $C$ is a projective curve aver $k$, its automorphism group
  $\aut C$ acts on the base $\M\gl$ of its universal deformation.
\end{lemma}
\begin{proof}
  Let $g \in \aut C$. As
  $$\english \xymatrix @!0 @R=1.5pc @C=3pc {
    && C \ar[rrdd]^\lambda \ar[ld]_{g^{-1}} \ar[ddd] \\
    & C \ar[ld]_{\lambda} \ar[rdd] \\
    \C \ar[ddd] &&&& \C \ar[ddd] \\
    && \spec k \ar[rrdd] \ar[lldd] \\
    \\
    \M &&&& \M} $$
  the universal property defining $\M\gl$ 
  gives us a cartesian square
  $$\english  \xymatrix{\relax
    \C \ar@{.>}[r] \ar[d] \cartesien & \C \ar[d] \\ 
    \M \ar@{.>}[r]_{\wt g} & \M } $$
  such that 
  $$\english \xymatrix @!0 @R=1.5pc @C=3pc {
    && C \ar[rrdd]^\lambda \ar[ld]_{g^{-1}} 
    \ar[ddd]|!{[dldl];[drdr]}\hole \\
    & C \ar[ld]_{\lambda} \\
    \C \ar[ddd] \ar@{.>}[rrrr] &&&& \C \ar[ddd] \\
    && \spec k \ar[rrdd] \ar[lldd] \\
    \\
    \M \ar@{.>}[rrrr]_{\wt g} &&&& \M.} $$
  We leave it to the reader to check that this actually defines a
  group action: use the universal property of the deformation 
  to show that 
  $\wt{gh}=\wt g \wt h$.
\end{proof}

\begin{lemma}\label{lemme:deformations=espacetangent}%
  Let $C$ be a stable $n$-pointed curve and 
  $\M\gl = \spec M\gl$ the base of its universal deformation.
  There is an $\aut C$-equivariant isomorphism
  between the tangent space to $\Mgnb$ at $[C]$ and
  that of $\M\gl$ at the origin, 
  $$ T_{[C]} \Mbgn \iso T_0 \M\gl. $$
\end{lemma}
\begin{proof}
  This results from the following isomorphisms of $\aut C$-equivariant
  vector spaces.
  \begin{align*}
    T_{[C]} \Mbgn &= 
    \cat{hom}_{[C]} \left(
      \spec\dfrac{k[\epsilon]}{\langle\epsilon^2\rangle}, 
      \Mbgn \right)
    \tag{1} \\
    &\iso D\left(\dfrac{k[\epsilon]}{\langle\epsilon^2\rangle}\right)
    \tag{2} \\
    &\iso \hom_k \left( M\gl,
      \dfrac{k[\epsilon]}{\langle\epsilon^2\rangle} \right). 
    \tag{3} \\
    &\iso \hom_k \left( 
      \spec \dfrac{k[\epsilon]}{\langle\epsilon^2\rangle},
      \M\gl
    \right)
    \tag{4} \\
    &= T_0 \M\gl.
    \tag{5}    
  \end{align*}
  The only delicate point is the isomorphism between (1) and (2): 
  instead of considering them as equivariant vector spaces, we
  shall regard them as categories and show that they are isomorphic. 

  The objects of the category 
  $\cat{hom}_{[C]}(
  \spec k[\epsilon]/\langle\epsilon^2\rangle,
  \Mbgn) $
  are the commutative diagrams
  $$\english \xymatrix @!0 @R=5pc @C=3pc {
    & \spec k \ar[ld] \ar[rd] ^{[C]} \\
    \spec \dfrac{k[\epsilon]}{\langle\epsilon^2\rangle}
    \ar[rr] && \Mgnb } $$
  and its morphisms are the 2-morphisms 
  $$
  \english\xymatrix{
    **[l] 
    \hbox to 2mm{\hss\ensuremath{
        \spec \dfrac{k[\epsilon]}{\langle\epsilon^2\rangle}}}
    \ar@/^2pc/[r]|*{}="A"
    \ar@/_2pc/[r]|*{}="B"
    \ar@{=>}"A";"B" 
    & \Mgnb. 
    } $$
  
  The objects of the category
  $D(k[\epsilon]/\langle\epsilon^2\rangle)$ 
  are classes of families of stable $n$-pointed genus $g$ curves over 
  $\spec k[\epsilon]/\langle\epsilon^2\rangle$ whose fiber at $0$
  is~$C$ : one identifies two such families if they differ by an 
  isomorphism whose restriction to $C$ is the identity; 
  its morphisms are isomorphisms of stable pointed curves. 

  Let us show that the functor 
  $$\english  \definefunction%
  {\cat{hom}_{[C]}\left(
      \spec \dfrac{k[\epsilon]}{\langle\epsilon^2\rangle},
      \Mbgn \right)}%
  {D\left( \dfrac{k[\epsilon]}{\langle\epsilon^2\rangle} \right)}%
  {\raisebox{.5\depth}{\xymatrix @!0 @R=5pc @C=3pc {
        & \spec k \ar[ld] \ar[rd]^{[C]} \\
        \spec \dfrac{k[\epsilon]}{\langle\epsilon^2\rangle}
        \ar[rr] && \Mbgn}}}%
  {\raisebox{.5\depth}{\xymatrix{C \ar[r]\ar[d]\cartesien & 
      \C \ar[d]\ar[r]\cartesien &
      \Mb_{g,n+1} \ar[d] \\
      \spec k \ar[r] & 
      \spec \dfrac{k[\epsilon]}{\langle\epsilon^2\rangle} \ar[r] &
      \Mbgn}}}
  $$
  is essentially surjective.
  Let 
  $$ \xymatrix{ C \ar[r] \ar[d] \cartesien & \C \ar[d] \\
    \spec k \ar[r] & 
    \spec \dfrac{k[\epsilon]}{\langle\epsilon^2\rangle} } $$
  be an object of $D(k[\epsilon]/\langle\epsilon^2\rangle)$. 
  The family of curves
  $\C \lra \spec k[\epsilon]/\langle\epsilon^2\rangle$
  is an object of the fiber 
  $\Mgnb(\spec k[\epsilon]/\langle\epsilon^2\rangle)$
  of the fibered category $\Mgnb$
  above 
  $\spec k[\epsilon]/\langle\epsilon^2\rangle$, 
  but this fiber is equivalent to 
  $$\cat{hom}(\spec k[\epsilon]/\langle\epsilon^2\rangle,
  \Mgnb). $$
  The composition 
  $ \spec k \lra \spec k[\epsilon]/\langle\epsilon^2\rangle \lra
  \Mgnb$ is given by the fiber product 
  $$\xymatrix{C \ar[r]\ar[d]\cartesien & \C \ar[d] \\
    \spec k \ar[r] & 
    \spec \dfrac{k[\epsilon]}{\langle\epsilon^2\rangle},
    } $$
  so it is (isomorphic to) $[C] : \spec k \lra \Mgnb$; 
  hence the family of curves 
  $\C \lra \spec k[\epsilon]/\langle\epsilon^2\rangle$ is 
  (isomorphic to) an object of
  $$\cat{hom}_{[C]} ( 
  \spec k[\epsilon]/\langle\epsilon^2\rangle, 
  \Mbgn).$$

  Let us now show that this functor is fully faithful. As the
  categories are groupoids, it suffices to prove that it induces
  isomorphisms between automorphism groups. But by definition, the
  automorphisms of 
  $$ \english \xymatrix @!0 @R=5pc @C=3pc { 
    & \spec k \ar[ld] \ar[rd]^{[C]} \\
    \spec \dfrac{k[\epsilon]}{\langle\epsilon^2\rangle}
    \ar[rr] && \Mgnb } $$
  are the automorphisms of 
  $$ \spec \dfrac{k[\epsilon]}{\langle\epsilon^2\rangle} \lra \Mgnb,
  $$
  \ie, (from Yoneda's lemma for fibered categories)
  the automorphisms of 
  $$\C \lra \spec \dfrac{k[\epsilon]}{\langle\epsilon^2\rangle}.$$

  To show that this equivalence of categories is actually an
  isomorphism, it suffices to remark that in both cases the set of morphisms
  from a given object is in bijection with the \emph{finite}
  set $\aut C$. 
\end{proof}

\begin{defi}
  A \d{graph} is a sextuplet
  $G = (F, V, j, \partial, g, c)$ where
  $F$ is a set (whose elements will be called \d{half-edges}), 
  $V$ is a set (whose elements will be called \d{vertices}),
  $j : F \lra F$ is an involution, 
  $\partial : F \lra V$ is a map.
  An \d{edge} of $G$ is a 2-element part of $F$ of the form 
  $\{\phi, j\phi\}$.
  The \d{legs} of $G$ are the fixed points of $j$. 
  Finally, $g: V \lra \NN$ is a map and 
  $c : \{\,\text{legs}\,\} \lra \NN$ is a bijection 
  between the set of legs and 
  $[\![1,n]\!]$ for some $n \in \NN$.

  A \d{graph morphism} 
  $(F,V,j,\partial,g,c) \lra (F',V',j',\partial',g',c')$ is the datum of two
  bijections 
  $f_1 : F \lra F'$ and $f_0 : V \lra V'$ such that 
  $f_1 j = j' f_1$, 
  $f_0 \partial = \partial' f_1$,
  $g=g'f_0$ and 
  $c = c' f_1$.

  The \d{graph associated} to a stable $n$-pointed curve $C$ (over an
  algebraically closed field) is
  \begin{align*}
    V &= \{\,\text{connected components of the normalization $\wt C$}\,\} \\
    &\iso \{\,\text{irreducible components of $C$}\,\} 
    \displaybreak[0]\\
    F &= \{\,\text{special points of $C$}\,\}
    \displaybreak[0]\\
    j&:{} \definefunction FF\phi{
      \hbox{\ensuremath{\phi'}\ \vtop{%
        \hsize=9cm\noindent\raggedright
        where
        $\pi^{-1}\pi\phi=\{\,\phi,\phi'\,\}$, 
        \mathstrut
        (one may have $\phi = \phi'$ if
        $\pi^{-1}\pi\phi$ has a single element, \ie, if $\phi$ is a
        \mathstrut
        special point above
        \mathstrut
        a marked point)}}
      }
    \displaybreak[0]\\    
    \partial &: \definefunction FV\phi{\text{connected component of $\wt
        C$ on which $\phi$ lies}}
    \displaybreak[0]\\
    g&:\definefunction V{\NN}v{\text{geometric genus of $v$}} \\
    c&:\definefunction{\{\, \text{legs} \,\}}{[\![1,n]\!]}{\phi}{i
      \text{ where $P_i = \pi \phi$ and the 
        $(P_i)_{1\leq i \leq n}$ are the marked points}}
  \end{align*}

  A \d{ribbon graph} is a graph endowed, at each vertex, with a cyclic
  order on the three outgoing half-edges. 
\end{defi}
  
  \begin{figure}[htbp]
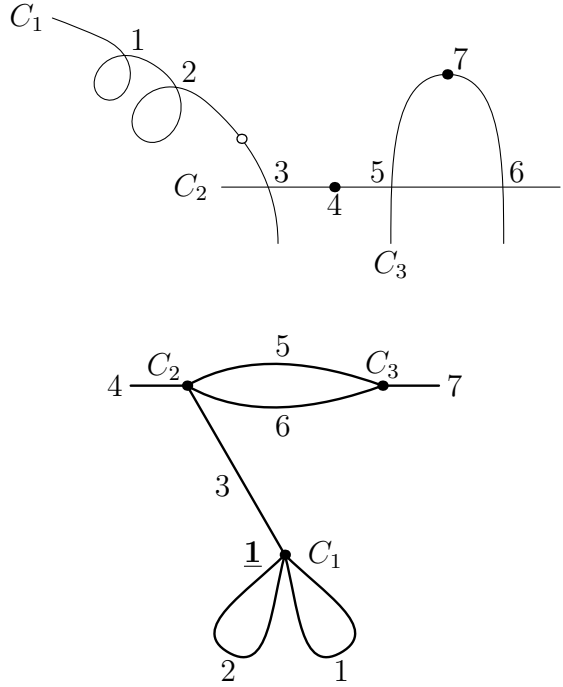

    \begin{center}
      \leavevmode
      \subfigure{\includegraphics{1.109}}%
      \qquad
      \subfigure{\includegraphics{1.110}}%
      \caption{A stable curve and the associated graph (the circle on
        $C_1$ indicates that its geometric genus is one; the genus of
        a vertex is indicated by a bold underlined number, if
        non-zero)} 
      \label{fig:correspondancecourbegraphe}
    \end{center}
  \end{figure}

\begin{theo}\label{theo:graphe-point-basetangenticiel}%
  Let $C$ be a maximally degenerate stable pointed curve
  corresponding to a point $[C]$ of $\Mgnb$ and let $G$ be the
  associated graph. Any ribbon graph structure on $G$ defines a
  tangencial base point on $\Mgn$ at $[C]$. 
\end{theo}
\begin{proof}
  Let $\wt G$ be a ribbon graph structure on $G$ and let 
  $x_1,\dots,x_n$ be the singular points of $C$. 
  
  This defines local coordinates in the neighborhood of the
  special points of $C$ in the following way. 
  A special point $x$ is a point of the normalization $\wt C$ 
  above a singular or a marked point of $C$ : it lies on a connected
  component of $\wt C$, that has exactly two other special points. 
  This component is a vertex of the graph $G$ and the three points are
  the three half-edges stemming from it.
  The ribbon graph structure gives us a cyclic order on these three
  points: one may define an isomorphism between the connected component of
  $\wt C$ and the projective line $\PP^1$, by sending the point $x$ to
  $0$, 
  the next one to $+1$ and the last one
  to $-1$; hence, we get coordinates around $x$. 
  These coordinates yield isomorphisms 
  $$ \wh \O_{C,\, x_i} \iso 
  \dfrac{k[\![u_i,v_i]\!]}{\langle u_i v_i \rangle}. $$
  
  As the universal deformation of 
  $k[\![u_i,v_i]\!] / \langle u_i v_i\rangle$ is
  $$\xymatrix{\relax
    \spec \dfrac{k[\![u_i,v_i]\!]}{\langle u_i v_i\rangle}
    \ar[r]\ar[d]\cartesien &
    \spec \dfrac{k[\![u_i,v_i,\epsilon_i]\!]}{\langle u_iv_i-\epsilon_i\rangle}
    \ar[d] \\
    \spec k \ar[r] & \spec k[\![\epsilon_i]\!], } $$
  the preceeding isomorphisms define an isomorphism
  $$ M\loc \iso k[\![\epsilon_1]\!] \hattens \cdots \hattens
  k[\![\epsilon_n]\!] = k[\![\epsilon_1,\dots,\epsilon_n]\!], $$
  hence, from lemma \ref{lemme:defgl=defloc},
  $M\gl \iso k[\![\epsilon_1,\dots,\epsilon_n]\!]. $
  
  The isomorphism 
  $ T_{[C]} \Mgnb \iso T_0 \M\loc $
  identifies the trace of the divisor $\Sgn$
  on the tangent space with the union of the hyperplanes 
  $\epsilon_i = 0$; from example \ref{ex:tangentspacetangbasepoint},
  the $\epsilon_i$ then define a tangencial base point. 
\end{proof}

\begin{rem}
  These tangencial base points were constructed in \cite{IN} in the
  case of a maximally degenerate pointed curve with no irreducible
  components of arithmetic genus 1.
\end{rem}

We shall conclude with a combinatorial description of this tangent
space. 

\begin{lemma}
  Let $C$ be s maximally degenerate stable $n$-pointed curve and $G$
  the associated graph. Then $\aut C \iso \aut G$. 
\end{lemma}
\begin{proof}
  Left to the reader: an automorphism $f$ of $C$ induces an automorphism
  $\wt f$ of the normalization $\wt C$ and hence an automorphism of
  the graph associated to $C$; check that this defines a group
  morphism $\aut C \lra \aut G$ and construct its inverse.
\end{proof}

\begin{coro}
  Let $C$ be a maximally degenerate stable pointed curve and $G$ the
  associated graph. A ribbon graph structure on $G$ defines a vector
  space isomorphism 
  $$ T_{[C]}\Mgnb
  \iso 
  V = \vect \{\, \textrm{\rm edges $e$ de $G$} \,\}. $$
  The action of $\aut C = \aut G$ 
  on this vector space is the following.
  If $g \in \aut G$ is an element of the group, $e$ an edge
  and $v_e$ the corresponding vector, we set 
  $$ g \cdot v_e = \pm v_{g \cdot e}. $$
  The sign is the following:
  \begin{itemize}
  \item[$+1$] if the automorphism $g$ preserves the cyclic order at
    both ends of $e$ or reverses it at both ends;
  \item[$-1$] otherwise. 
  \end{itemize}
\end{coro}
\begin{proof}
  Left to the reader: take local coordinates as in the proof of 
  \ref{theo:graphe-point-basetangenticiel} and use them to describe
  the action of $\aut C$. 
\end{proof}

\cbend

\bibliographystyle{amsplain}
\nocite{*}
\bibliography{1}
\end{document}